\documentclass[fleqn]{mat01}
\usepackage{times,mathtimy,amssymb}
\begin{document}

\font\xxx=tir at 10pt
\def\d{\hbox{\xxx{d}}}

\font\zzz=mtmib at 13.5pt
\def\N{\hbox{\zzz{N}}}
\def\r{\hbox{\zzz{r}}}

\font\zzzzz=mtmib at 9.4pt
\def\q{\hbox{\zzzzz{q}}}

\newtheorem{definit}{\rm DEFINITION}
\newtheorem{theo}{Theorem}
\renewcommand{\thetheo}{\arabic{theo}}
\newtheorem{theor}[theo]{\bf Theorem}
\newtheorem{corolla}{\rm COROLLARY}
\newtheorem{lem}{Lemma}

\setcounter{page}{223}
\firstpage{223}

\title{On the absolute $\N_{\q_{\alpha}}$-summability of $\r$th derived
conjugate series}

\markboth{A~K~Sahoo}{On the absolute $N_{q_{\alpha}}$-summability of $r$th derived
conjugate series}

\author{A~K~SAHOO}

\address{Department of Mathematics, Government Kolasib College,
P.~B.~No.~23, Kolasib~796~081, India}

\volume{113}
\renewcommand\thesubsection{\arabic{section}.\arabic{subsection}.}
\mon{August}

\parts{3}

\Date{MS received 26 March 2002}

\begin{abstract}
The object of the present paper is to study the absolute
$N_{q_{\alpha}}$-summability of $r$th derived conjugate series
generalizing a known result.
\end{abstract}

\keyword{Fourier series; conjugate series; derived conjugate series;
Nevanlinna summability; kernel.}

\maketitle

\section{Introduction}

\subsection{}

In the year 1921, Nevanlinna \cite{7} suggested and discussed an interesting
method of summation called $N_{q}$\!-method. Moursund~\cite{5} applied this
method for summation of Fourier series and its conjugate series. Later,
Moursund~\cite{6} developed $N_{q_{p}}$\!-method (where $p$ is a positive
integer) and applied it to $p$th derived Fourier series. Samal~\cite{9}
discussed $N_{q_{\alpha}}$\!-method $(0 \leq \alpha < 1)$ and studied
absolute $N_{q_{\alpha}}$\!-summability of some series associated with
Fourier series. In his Ph.D. thesis \cite{10} he extended $N_{q_{p}}$\!-method
of summation to $N_{q_{\alpha}}$\!-method for any $\alpha \geq 0$ and
studied absolute $N_{q_{\alpha}}$\!-summability of Fourier series. Earlier
we \cite{8} have studied absolute $N_{q_{\alpha}}$\!-summability of a series
conjugate to a Fourier series. In the present paper we shall study the
absolute $N_{q_{\alpha}}$\!-summability of $r$th $(r < \alpha)$ derived
series of a conjugate series.

\subsection{}

\begin{definit}{\rm \cite{6,10}}$\left.\right.$\vspace{.5pc}

\noindent {\rm Let $F(w)$ be a function of a continuous parameter $w$ defined for all
$w > 0$. The $N_{q_{\alpha}}$\!-method consists in forming the
$N_{q_{\alpha}}$\!-transform or mean
\begin{equation*}
N_{q_{\alpha}} F(w) = \int_{0}^{1} q_{\alpha} (t) F(wt) \d t
\end{equation*}
and then considering the limit
\begin{equation*}
\lim\limits_{w \rightarrow \infty} N_{q_{\alpha}} F(w),
\end{equation*}

\pagebreak

\noindent where the class of functions $q_{\alpha} (t)$ is such that}\vspace{-.4pc}
\end{definit}

\begin{enumerate}
\renewcommand{\labelenumi}{(\arabic{enumi})}
\leftskip .1pc
\item $q_{\alpha} (t) \geq 0$ for $0 \leq t \leq 1$,\vspace{.4pc}

\item $\displaystyle \int_{0}^{1} q_{\alpha} (t) \d t = 1$,\vspace{.6pc}

\item $\displaystyle \frac{\d^{\beta}}{\d t^{\beta}}
q_{\alpha} (t)$ exists and is absolutely continuous for $0 \leq t \leq
1, \beta = 1, 2, \ldots, k - 1$, where $[\alpha] = k$,\vspace{.4pc}

\item $\displaystyle \frac{\d^{\beta}}{\d t^{\beta}}
q_{\alpha} (t) = 0$ for $t = 1, \beta = 0, 1, 2, \ldots, k - 1$,\vspace{.4pc}

\item $\displaystyle \frac{\d^{k}}{\d t^{k}} q_{\alpha} (t)$ exists for
$0 < t < 1$,\vspace{.4pc}

\item $\displaystyle (-1)^{k} \frac{\d^{k}}{\d t^{k}} q_{\alpha} (t)
\geq 0$ and monotonic increasing for $0 < t < 1$,\vspace{.6pc}

\item $\displaystyle \int_{0}^{t} \frac{Q_{k} (u)}{u^{1 + \alpha -
k}} \d u = O \left( \frac{Q_{k} (t)}{t^{\alpha - k}} \right)$,\vspace{.4pc}

where
\begin{equation*}
\hskip -1.2pc Q_{k} (t) = \int_{1-t}^{1} (-1)^{k} \frac{\d^{k}}{\d u^{k}}
q_{\alpha} (u)\d u.
\end{equation*}
\end{enumerate}
Also we set
\begin{equation*}
Q(t) = \int_{1-t}^{1} q_{\alpha} (u) \d u.
\end{equation*}
If $\lim_{w \rightarrow \infty} N_{q_{\alpha}} F(w)$ exists, we
say that $N_{q_{\alpha}}$\!-limit of $F(w)$ exists.

\begin{definit}{\rm \cite{9,10}}$\left.\right.$\vspace{.5pc}

\noindent {\rm Let $\sum_{n = 0}^{\infty} u_{n}$ be an infinite series
with $S(w) = \sum_{n \leq w} u_{n}$. If $\lim_{w \rightarrow \infty}
\{ \sum_{n \leq w} u_{n} Q(1 - (n/w))\} = 1$, we say
that $\sum u_{n}$ is summable by $N_{q_{\alpha}}$\!-method to the sum~1.
In short we write that $\sum u_{n} = 1 (N_{q_{\alpha}})$. Further the
series $\sum u_{n}$ is said to be $|N_{q_{\alpha}} |$-summable (absolute
$N_{q_{\alpha}}$\!-summable) if
\begin{equation*}
\int_{A}^{\infty} \frac{\d w}{w^{2}} \left\vert \sum\limits_{n
\leq w} nu_{n} q_{\alpha} \left( \frac{n}{w} \right) \right\vert <
\infty
\end{equation*}
for some positive constant $A$.

For $\alpha = 0$, the method reduces to the original $N_{q}$\!-method \cite{7}
and if $\alpha$ is any positive integer $p$, then the method reduces to
$N_{q_{p}}$\!-method of Moursund \cite{6}.}
\end{definit}

\subsection{}

Let $f(t)$ be a periodic function with period $2\pi$ and Lebesque
integrable over $(-\pi, \pi)$.

Let
\renewcommand{\theequation}{\thesubsection\arabic{equation}}
\setcounter{equation}{0}
\begin{equation}
f(t) \sim \frac{1}{2} a_{0} + \sum\limits_{n = 1}^{\infty} (a_{n}\ 
\cos\ nt + b_{n}\ \sin\ nt) \equiv \sum\limits_{n = 0}^{\infty} A_{n}
(t).
\end{equation}
The series conjugate to (1.3.1) at $t = x$ is given by
\begin{align}
&\sum\limits_{n = 1}^{\infty} (b_{n}\ \cos\ nx - a_{n}\ \sin\ nx) \equiv
\sum\limits_{n = 1}^{\infty} B_{n}(x),\\
&P(u) = \sum\limits_{i = 0}^{r - 1} \frac{\theta_{i}}{i!} u^{i}\quad
\hbox{for}\quad -\pi \leq u \leq \pi,\nonumber
\end{align}
where $\theta_{i}$s for $i = 0, 1, 2, \ldots, r-1$ are arbitrary constants.
\begin{align*}
&h(u) = \frac{\{ f(x + u) - P(u)\} - (-1)^{r} \{f(x - u) - P(-
u)\}}{2u^{r}},\\
&H_{0} (t) = h(t),\\
&H_{\beta} (t) = \frac{1}{\Gamma (\beta)} \int_{0}^{t} (t -
u)^{\beta - 1} h(u) \d u,\quad (\beta > 0),\\
&h_{\beta} (t) = \Gamma (1 + \beta) t^{-\beta} H_{\beta} (t),\quad
(\beta \geq 0).
\end{align*}

\section{Purpose of the present paper}

In the present paper we shall prove the following theorems:

\begin{theor}[\!]
Let $\beta = \alpha - r$. If $H_{\beta} (+ 0) = 0$ and
$\int_{0}^{\pi} t^{-\beta} |\d H_{\beta} (t)| < \infty${\rm ,} where $\beta >
0${\rm ,} then the $r$th derived series of the conjugate series of $f(t)$ at
$t = x$ is $|N_{q_{\alpha}}|$\!-summable.
\end{theor}

\begin{theor}[\!]
Let $\rho = \alpha - r - 1$. If $\rho \geq 0$ and $\int_{0}^{\pi} t^{-1}
|h_{\rho} (t)| \d t < \infty${\rm ,} then the $r$th derived series of the
conjugate series of $f(t)$ at $t = x$ is $|N_{q_{\alpha}}|$-summable.
\end{theor}

By taking $\beta = \rho + 1, \rho \geq 0$ in Theorem~1, we can obtain
Theorem~2 at once as it is known \cite{4} that
\begin{align*}
&H_{\rho + 1} (+0) = 0\quad\hbox{and}\quad \int_{0}^{\pi} t^{-\rho
- 1} |\d H_{\rho + 1} (t)| < \infty\\
&\quad\ \Longleftrightarrow h_{\rho + 1} (t) \in BV (0, \pi)\quad
\hbox{and}\quad \frac{h_{\rho + 1} (t)}{t} \in L(0, \pi)\\
&\quad\ \Longleftrightarrow \frac{h_{\rho} (t)}{t} \in L(0, \pi).
\end{align*}

By taking $q_{\alpha} (t) = (\alpha + \delta) (1 - t)^{\alpha + \delta -
1}$, where $\delta > 0$ and $\alpha + \delta < k + 1\ ([\alpha] = k)$ in
Theorems~1 and 2, we obtain the following corollaries respectively.

\begin{corolla}{\rm \cite{3}}$\left.\right.$\vspace{.5pc}

\noindent If $H_{\beta} (+0) = 0$ and $\int_{0}^{\pi} 
t^{-\beta} |\d H_{\beta} (t)| < \infty${\rm ,} then the $r$th derived series
of the conjugate series of $f(t)$ at $t = x$ is summable $|C, \beta + r
+ \delta|${\rm ,} where $\beta > 0$ and $\delta > 0$.
\end{corolla}

\begin{corolla}{\rm \cite{3}}$\left.\right.$\vspace{.5pc}

\noindent If $\rho \geq 0$ and $\int_{0}^{\pi} t^{-1} |h_{\rho} (t)| \d t < \infty${\rm ,} then the 
$r$th derived series of the conjugate series of $f(t)$ at $t = x$ is
summable $|C, \rho + r + 1 + \delta|$.
\end{corolla}

\section{Notations and lemmas}

\subsection{\it Notations}

For our purpose we use the following notations throughout this paper.
\begin{align*}
&[\alpha] = k,\\
&m = \min\ (k - r, r),\\
&q^{k} (u) = (-1)^{k} \frac{\d^{k}}{\d u^{k}}\ q_{\alpha} (u),\\
&(\cos\ nu)_{j} = \left( \frac{\d}{\d u} \right)^{j}\ \cos\ nu,\\
&S^{i,j} (x,u) = \sum\limits_{n \leq x} (x - n)^{i} (\cos\ nu)_{j},\\
&G_{i} (w, u) = \sum\limits_{n \leq w} q_{\alpha} \left( \frac{n}{w}
\right) \left( \frac{\d}{\d u}\right)^{k + 1 - i}\ \cos\ nu, \quad
\hbox{for}\quad i = 0, 1, 2, \ldots, m,\\
&g_{i} (x, w, u) = \frac{1}{k!} (-1)^{k} \left( \frac{\d}{\d x}
\right)^{k} q_{\alpha} \left( \frac{x}{w} \right) \frac{\d}{\d x}\ S^{k,
k+1-i} (x, u)\\
&\hskip 2.5cm \hbox{for}\quad i = 0, 1, 2, \ldots, m.
\end{align*}

\subsection{\it Lemmas}

We need the following lemmas for the proof of our theorem.

\def\lemaa{\trivlist\item[\hskip\labelsep{\it Lemma {\rm 1 [1]}.}]}
\begin{lemaa}{\it 
If $\beta > \alpha > 0, H_{\alpha} (t)$ is of $BV(0, \pi)$ and
$H_{\alpha} (+0) = 0${\rm ,} then $H_{\beta} (t)$ is an integral in $(0, \pi)$
and for almost all values of $t${\rm ,}
\begin{equation*}
H_{\beta}' (t) = \frac{1}{\Gamma (\beta - \alpha)} \int_{0}^{t}
(t - u)^{\beta - \alpha - 1} \d H_{\alpha} (u).
\end{equation*}}
\end{lemaa}\vspace{-.3pc}

\def\lemaa{\trivlist\item[\hskip\labelsep{\it Lemma {\rm 2 [6]}.}]}
\begin{lemaa}{\it
If $\alpha \geq 1${\rm ,} the kernel $q_{\alpha} (t)$ is monotonic decreasing{\rm ,}
its derivatives of odd orders less than $k$ are negative and monotonic
increasing{\rm ,} its derivatives of even orders less than $k$ are positive
and monotonic decreasing and there exists a constant $A_{k}$ such that
\begin{equation*}
\left\vert \frac{\d^{\beta}}{\d t^{\beta}}\ q_{\alpha} (t) \right\vert <
A_{k}\quad (\beta = 0, 1, 2, \ldots, k-1)
\end{equation*}
and
\begin{equation*}
\int_{0}^{1} \left\vert \frac{\d^{k}}{\d t^{k}}\ q_{\alpha} (t)
\right\vert \d t < A_{k}.
\end{equation*}}
\end{lemaa}\vspace{-.3pc}

\def\lemaa{\trivlist\item[\hskip\labelsep{\it Lemma {\rm 3 [10]}.}]}
\begin{lemaa}{\it 
$Q_{k} (t)$ is continuous and monotonic increasing function of $t${\rm ,}
$Q_{k} (t) \geq 0, Q(0) = 0$ and $Q (1) = 1$.}
\end{lemaa}\vspace{.5pc}

This follows directly from the definition of $Q(t)$ and $Q_{k}(t)$.

\def\lemaa{\trivlist\item[\hskip\labelsep{\it Lemma {\rm 4 [10,8]}.}]}
\begin{lemaa}{\it $\int_{0}^{1} q^{k}(t)/((1-t)^{\alpha -
k}) \d t$ exists.}
\end{lemaa}

\pagebreak

\def\lemaa{\trivlist\item[\hskip\labelsep{\it Lemma {\rm 5 [8]}.}]}
\begin{lemaa}{\it 
Let $x > 0$.}\vspace{-.4pc}
\end{lemaa}

{\it \begin{enumerate}
\renewcommand{\labelenumi}{\rm (\roman{enumi})}
\leftskip .1pc
\item If $1/x < u \leq\pi${\rm ,} then
\begin{equation*}
\hskip -1.5pc S^{i,j} (x, u) = \begin{cases} 
O(x^{i} u^{-j-1}) &\hbox{for}\quad 0 \leq j \leq i,\\[.2pc]
O(x^{j} u^{-i-1}) &\hbox{for}\quad j > j \geq 0.
\end{cases}
\end{equation*}

\item If $1/x \geq u > 0${\rm ,} then
\begin{equation*}
\hskip -1.5pc S^{i,j} (x, u) = O(x^{i+j+1}).
\end{equation*}
\end{enumerate}}

\def\lemaa{\trivlist\item[\hskip\labelsep{\it Lemma {\rm 6 [2]}.}]}
\begin{lemaa}{\it 
Let $\lambda = \{\lambda_{n}\}$ be a positive monotonic increasing
sequence with $\lambda_{n} \rightarrow \infty$ as $n \rightarrow
\infty$. Then
\begin{equation*}
A_{\lambda} (x) = A_{\lambda}^{0} (x) = \sum\limits_{\lambda_{n} \leq x}
a_{n}
\end{equation*}
and
\begin{equation*}
A_{\lambda}^{r} (x) = \sum\limits_{\lambda_{n} \leq x} (x -
\lambda_{n})^{r} a_{n} (r > 0).
\end{equation*}
If $k$ is a positive integer{\rm ,}
\begin{equation*}
A_{\lambda} (x) = \frac{1}{k!} \left( \frac{\d}{\d x} \right)^{k}
A_{\lambda}^{k} (x).
\end{equation*}}
\end{lemaa}

\def\lemaa{\trivlist\item[\hskip\labelsep{\it Lemma {\rm 7 [8,10]}.}]}
\begin{lemaa}{\it
For $\alpha \geq 1${\rm ,}
\begin{equation*}
\sum\limits_{n \leq w} (-1)^{n} n^{k} q_{\alpha} \left( \frac{n}{w}
\right) = O \left\lbrace q^{k} \left( 1 -\frac{1}{w} \right)
\right\rbrace + O \left\lbrace wQ_{k} \left( \frac{1}{w} \right)
\right\rbrace.
\end{equation*}}
\end{lemaa}

\def\lemaa{\trivlist\item[\hskip\labelsep{\it Lemma {\rm 8 [8,10]}.}]}
\begin{lemaa}{\it
For $\alpha \geq 1$ and $r =0, 1, 2, \ldots, k-1${\rm ,}
\begin{equation*}
\sum\limits_{n \leq w} (-1)^{n} n^{r} q_{\alpha} \left( \frac{n}{w}
\right) = O(1).
\end{equation*}}
\end{lemaa}

\setcounter{lem}{8}
\begin{lem}
For $i = 0, 1, 2, \ldots, k-1${\rm ,}
\begin{equation*}
\sum\limits_{n \leq w} (-1)^{n} n^{i} \in |N q_{\alpha}|.
\end{equation*}
\end{lem}

\begin{proof}
For $i = 0, 1, 2, \ldots, k-2$,
\begin{align*}
\int_{1}^{\infty} \frac{\d w}{w^{2}} \left\vert \sum\limits_{n
\leq w} n(-1)^{n} n^{i} q_{\alpha} \left( \frac{n}{w} \right)
\right\vert &= \int_{1}^{\infty} O(1) \frac{\d w}{w^{2}}\quad \hbox{by 
Lemma~8}\\
&= O(1)
\end{align*}
and
\begin{align*}
&\int_{1}^{\infty} \frac{\d w}{w^{2}} \left\vert \sum\limits_{n
\leq w} (-1)^{n} n^{k} q_{\alpha} \left( \frac{n}{w} \right)
\right\vert\\[.3pc]
&\quad\ = \int_{1}^{\infty} O\left\lbrace q^{k} \left( 1 - \frac{1}{w}
\right) \right\rbrace \frac{\d w}{w^{2}}\\[.3pc]
&\qquad\ \ + \int_{1}^{\infty} O
\left\lbrace wQ_{k} \left( \frac{1}{w} \right) \right\rbrace \frac{\d
w}{w^{2}}\quad \hbox{by Lemma~7}\\[.3pc]
&\quad\ = O\left( \int_{0}^{1} q^{k} (u) \d u \right) + O \left(
\int_{0}^{1} \frac{Q_{k} (u)}{u}\ \d u \right)\\[.3pc]
&\quad\ =O(1)
\end{align*}
by Lemma~2 and the definitions of $q^{k} (u)$ and $Q_{k} (u)$.
This completes the proof of Lemma~9.
\end{proof}

\begin{lem}
For $i = 0, 1, 2, \ldots, m${\rm ,}
\begin{equation*}
G_{i} (w, u) = \int_{1}^{w} g_{i} (x, w, u) \d x.
\end{equation*}
\end{lem}

\begin{proof}
\begin{align*}
&G_{i} (w, u) = \sum\limits_{n \leq w} q_{\alpha} \left( \frac{n}{w}
\right) \left( \frac{\d}{\d u} \right)^{k + 1 -i} \cos\,nu\\[.3pc]
&\quad\ = q_{\alpha} (1) \sum\limits_{n \leq w} \left( \frac{\d}{\d u}
\right)^{k + 1 - i}\cos\,nu\\[.3pc]
&\qquad\ \ - \int_{1}^{w} \frac{\d}{\d x}\ 
q_{\alpha} \left( \frac{x}{w} \right) \left\lbrace \sum\limits_{n \leq
x} \left( \frac{\d}{\d u} \right)^{k + 1 - i} \cos\,nu\right\rbrace \d
x\\[.3pc]
&\quad\ = -\int_{1}^{w} \frac{\d}{\d x}\ q_{\alpha} \left( \frac{x}{w}
\right) \frac{1}{k!} \left( \frac{\d}{\d x} \right)^{k} \left\lbrace
\sum\limits_{n \leq x} (x - n)^{k} \left( \frac{\d}{\d u} \right)^{k + 1
- i} \cos\,nu \right\rbrace \d x\\[.3pc]
&\qquad\ \ \hbox{by Lemma~6}\\[.3pc]
&\quad\ = \left\lbrack \frac{1}{k!} \sum\limits_{\rho = 1}^{k - 1} (-1)^{\rho}
\left( \frac{\d}{\d x} \right)^{\rho} q_{\alpha} \left( \frac{x}{w}
\right) \left( \frac{\d}{\d x} \right)^{k - \rho}\ S^{k, k+1-i}\ (x,u)
\right\rbrack^{w}_{x=1}\\[.3pc]
&\qquad\ \ + \int_{1}^{w} \frac{(-1)^{k}}{k!} \left( \frac{\d}{\d x}
\right)^{k} q_{\alpha} \left( \frac{x}{w} \right) \frac{\d}{\d x}\ S^{k,
k+1-i}\ (x,u) \d x\\[.3pc]
&\qquad\ \ \hbox{(integrating by parts for $(k-1)$ times)}\\[.3pc]
&\quad\ = \int_{1}^{w} g_{i} (x,w,u)\d x
\end{align*}
as the integrated part vanishes for $x = w$ and $x = 1$.
\end{proof}

\begin{lem}
For $wt \leq \pi$ and $i = 0, 1, 2, \ldots, m${\rm ,}
\begin{equation*}
\int_{t}^{t+(1/w)} u^{r-i} (u-t)^{k-\alpha} G_{i} (w, u) \d u = O
(w^{\alpha - r + 1}).
\end{equation*}\vspace{.3pc}
\end{lem}

\begin{proof}
For $i = 0, 1, 2,\ldots, m$,
\begin{align*}
&\int_{t}^{t+(1/w)} u^{r-i} (u-t)^{k -\alpha} G_{i} (w, u)\d u\\
&\quad\ = \int_{t}^{t+(1/w)} u^{r-i} (u-t)^{k -\alpha} \left(
\sum\limits_{n \leq w} q_{\alpha} \left( \frac{n}{w} \right) \left(
\frac{\d}{\d u} \right)^{k + 1-i} \cos\,nu \right) \d u\\
&\quad\ = \int_{t}^{t+(1/w)} u^{r-i} (u-t)^{k -\alpha} O(w^{k+2-i}) \d
u\\
&\quad\ = O \left\lbrace \left(t + \frac{1}{w} \right)^{r - i} w^{k+2-i}
\int_{t}^{t + (1/w)} (u - t)^{k-\alpha} \d u \right\rbrace\\
&\quad\ = O \left\lbrace \left( \frac{wt + 1}{w} \right)^{r-i} w^{k+2-i} \cdot
\frac{1}{w^{k-\alpha+1}} \right\rbrace\\
&\quad\ = O (w^{\alpha -r + 1}) \quad \hbox{as}\quad wt \leq \pi.
\end{align*}
\end{proof}

\begin{lem}
For $i = 0, 1, 2,\ldots, m$ and $wt \leq \pi${\rm ,}
\begin{equation*}
\int_{t + (1/w)}^{\pi} u^{r-i} (u - t)^{k -\alpha} G_{i} (w, u)
\d u = O(w^{\alpha - r + 1}).
\end{equation*}
\end{lem}\vspace{.2pc}

\begin{proof}
By the use of Lemma~10,
\setcounter{equation}{0}
\begin{align}
&\int_{t + (1/w)}^{\pi} u^{r-i} (u - t)^{k -\alpha} G_{i} (w, u)
\d u\nonumber\\
&\quad\ = \int_{t + (1/w)}^{\pi} u^{r-i} (u - t)^{k -\alpha} \d u
\int_{1}^{w} g_{i} (x, w, u) \d x\nonumber\\
&\quad\ = \int_{t + (1/w)}^{\pi} u^{r-i} (u - t)^{k -\alpha} \d u
\frac{1}{k!}\nonumber\\
&\qquad\ \ \times \int_{1}^{w} (-1)^{k} \left( \frac{\d}{\d x}
\right)^{k} q_{\alpha} \left( \frac{x}{w} \right) \frac{\d}{\d x} S^{k,
k+1-i} (x,u) \d x\nonumber\\
&\quad\ = \frac{1}{(k-1)!} \int_{1}^{w} (-1)^{k} \left( \frac{\d}{\d x}
\right)^{k} q_{\alpha} \left( \frac{x}{w} \right) \d x\nonumber\\
&\qquad\ \ \times \int_{t + (1/w)}^{\pi} u^{r-i} (u - t)^{k -\alpha}\ S^{k-1,k+1-i}
(x,u) \d u\nonumber\\
&\quad\ = \frac{1}{(k-1)!} \int_{1}^{w} (-1)^{k} \left( \frac{\d}{\d x}
\right)^{k} q_{\alpha} \left( \frac{x}{w} \right) \d x\ w^{\alpha -k}\nonumber\\
&\qquad\ \ \times \int_{t + (1/w)}^{\xi} u^{r-i}\ S^{k-1,k+1-i}\ (x,u) \d u,
\end{align}
for some $t + (1/w) < \xi < \pi$, by an application of the mean value
theorem. For $i \geq 2$, using Lemma~5(i) in (3.2.1), we get
\begin{align*}
&\int_{t + (1/w)}^{\pi} u^{r-i} (u-t)^{k-\alpha} G_{i} (w,u) \d
u\\[.2pc]
&\quad\ = \frac{1}{(k-1)!} \int_{1}^{w} (-1)^{k} \left(
\frac{\d}{\d x} \right)^{k} q_{\alpha} \left( \frac{x}{w} \right)
w^{\alpha - k} \d x\\[.2pc]
&\qquad\ \ \times \int_{t + (1/w)}^{\xi} u^{r-i} O(x^{k-1}
u^{-k-2+i})\d u\\[.2pc]
&\quad\ = \frac{1}{(k-1)!} \int_{1}^{w} (-1)^{k} \left(
\frac{\d}{\d x} \right)^{k} q_{\alpha} \left( \frac{x}{w} \right)
w^{\alpha - k} O \left\lbrace \frac{x^{k-1}}{\left( t + \frac{1}{w}
\right)^{k-r+1}} \right\rbrace \d x\\[.2pc]
&\quad\ = O\left( w^{\alpha - r + 1} \int_{1}^{w} (-1)^{k} \left(
\frac{\d}{\d x} \right)^{k} q_{\alpha} \left( \frac{x}{w} \right) x^{k-1}
\d x \right)\\[.2pc]
&\quad\ = O \left(w^{\alpha - r + 1} \int_{0}^{1} q^{k} (\theta)
\d\theta \right)\\[.2pc]
&\quad\ = O \left(w^{\alpha - r + 1}\right)\quad \hbox{by Lemma~2.}
\end{align*}
For $i = 1$,
\begin{align*}
&\int_{t + (1/w)}^{\xi} u^{r-i}\ S^{k-1,k+1-i}\ (x,u) \d u\\[.2pc]
&\quad\ = \int_{t + (1/w)}^{\xi} u^{r-1}\ S^{k-1,k}\ (x,u) \d u\\[.2pc]
&\quad\ = \left\lbrack u^{r-1}\ S^{k-1,k-1}\ (x,u)\right\rbrack_{t +
(1/w)}^{\xi}\\[.2pc]
&\qquad\ \ -(r\!-\!1) \int_{t+(1/w)}^{\xi}  u^{r-2}\ S^{k-1,k-1}\ (x,u)\d u\\[.2pc]
&\quad\ = O \left\lbrace \frac{x^{k-1}}{\left( t + \frac{1}{w}
\right)^{k-r+1}} \right\rbrace + \int_{t + (1/w)}^{\xi} u^{r - 2}
O (x^{k-1} u^{-k})\d u\\
&\hskip 8cm \hbox{by Lemma~5(i)}\\
&\quad\ = O (w^{2k-r}).
\end{align*}
Similarly, for $i =0$, integrating by parts twice and using Lemma~5(i),
it follows that
\begin{equation*}
\int_{t + (1/w)}^{\xi} u^{r-i} S^{k-1,k+1-i} (x, u) \d u =
O(w^{2k-r}).
\end{equation*}
Hence, for $i \leq 1$, using the above estimation in (3.2.1)
\begin{align*}
&\int_{t + (1/w)}^{\pi} u^{r-i} (u-t)^{k-\alpha} G_{i} (w, u) \d
u\\[.2pc]
&\quad\ = O\left( \int_{1}^{w} (-1)^{k} \left( \frac{\d}{\d x}
\right)^{k} q_{\alpha} \left( \frac{x}{w} \right) w^{\alpha + k - r} \d
x \right)\\[.2pc]
&\quad\ = O \left( w^{\alpha-r+1} \int_{0}^{1} q^{k} (\theta) \d
\theta \right)\\[.2pc]
&\quad\ = O (w^{\alpha - r+1})\quad \hbox{by Lemma~2.}
\end{align*}
This completes the proof of Lemma~12.
\end{proof}

\begin{lem}
\begin{align*}
&\int_{t}^{t + (1/w)} u^{r-i} (u - t)^{k-\alpha} \d u
\int_{1}^{w - (\pi/t)} g_{i} (x,w,u)\d x\\[.2pc]
&\quad\ = O \left(
\frac{w^{\alpha - k}}{t^{k+1-r}}\ q^{k} \left( 1-\frac{\pi}{wt} \right)
\right).
\end{align*}
\end{lem}

\begin{proof}
For some $1 < \xi < w - (\pi/t)$, by an application of the mean value
theorem,
\begin{align}
&\int_{t}^{t+(1/w)} u^{r-i}(u-t)^{k-\alpha} \d u
\int_{1}^{w-(\pi/t)} g_{i} (x,w,u)\d x\nonumber\\[.2pc]
&\quad\ = \int_{t}^{t+(1/w)} u^{r-i}(u-t)^{k-\alpha} \d u\nonumber\\[.2pc]
&\qquad\ \ \times \int_{1}^{w-(\pi/t)} \frac{(-1)^{k}}{k!} \left(\frac{\d}{\d
x}\right)^{k}q_{\alpha} \left(\frac{x}{w}\right)\frac{\d}{\d x}
S^{k,k+1-i} (x,u)\d x\nonumber\\[.2pc]
&\quad\ = \int_{t}^{t+(1/w)} \frac{1}{k!} u^{r-i}(u-t)^{k-\alpha}
\left[(-1)^{k}\left(\frac{\d}{\d x}\right)^{k} q_{\alpha}
\left(\frac{x}{w}\right)\right]_{x = w-(\pi/t)} \d u\nonumber\\[.2pc]
&\qquad\ \ \times \int_{\xi}^{w-(\pi/t)} \frac{\d}{\d x} S^{k,k+1-i} (x,u)\d x\nonumber\\[.2pc]
&\quad\ = \frac{1}{k!} \int_{t}^{t+(1/w)} \frac{u^{r-i}(u\!-\!t)^{k-
\alpha}}{w^{k}}\ q^{k} \left(1\!-\!\frac{\pi}{wt}\right) \left[S^{k,k+1-
i}(x,u)\right]_{x=\xi}^{w-(\pi/t)}\d u.
\end{align}
For $i = 0$, using Lemma~5(i) in (3.2.2), we get
\begin{align*}
&\int_{t}^{t+(1/w)} u^{r-i}(u-t)^{k-\alpha}\d u
\int_{1}^{w-(\pi/t)} g_{i} (x,w,u)\d x\\[.2pc]
&\quad\ = \frac{1}{k!} \int_{t}^{t+(1/w)} \frac{u^{r}(u-t)^{k-
\alpha}}{w^{k}}\ q^{k} \left(1-\frac{\pi}{wt}\right)
O\left\lbrace\frac{\left(w-
\frac{\pi}{t}\right)^{k+1}}{u^{k+1}}\right\rbrace \d u
\end{align*}
\begin{align*}
&\quad\ = O\left(\frac{wq^{k}\left(1-\frac{\pi}{wt}\right)}{t^{k+1-r}}
\int_{t}^{t+(1/w)}(u-t)^{k-\alpha}\ \d u \right)\\[.2pc]
&\quad\ = O\left(\frac{w^{\alpha-k}q^{k}\left(1-\frac{\pi}{kt}\right)}{t^{k+1-
r}}\right).
\end{align*}
For $i \geq 1$, using Lemma~5(i) in (3.2.2), we obtain
\begin{align*}
&\int_{t}^{t+(1/w)} u^{r-i}(u-t)^{k-\alpha}\d u
\int_{1}^{w-(\pi/t)} g_{i} (x,w,u)\d x\\[.2pc]
&\quad\ = \frac{1}{k!} \int_{t}^{t+(1/w)} \frac{u^{r-i}(u-t)^{k-
\alpha}}{w^{k}} q^{k} \left(1-\frac{\pi}{wt}\right)
O\left(\frac{w^{k}}{u^{k+2-i}}\right) \d u\\[.2pc]
&\quad\ = O\left(\frac{q^{k}\left(1-\frac{\pi}{wt}\right)}{t^{k+2-r}}
\int_{t}^{t+(1/w)}(u-t)^{k-\alpha} \d u \right)\\[.2pc]
&\quad\ = O\left(\frac{w^{\alpha-k+1}}{t^{k+2-r}} q^{k}\left(1-
\frac{\pi}{wt}\right)\right).\\[.2pc]
&\quad\ = O\left(\frac{w^{\alpha-k}}{t^{k+1-r}} q^{k}\left(1-
\frac{\pi}{wt}\right)\right)\quad \hbox{as}\quad wt > \pi.
\end{align*}
This completes the proof of Lemma~13.\vspace{.4pc}
\end{proof}

\begin{lem}
For $i = 0,1,2,\ldots,m${\rm ,}
\begin{equation*}
\int_{t+(1/w)}^{\pi} u^{r-i} (u-t)^{k-\alpha}\ S^{k,k+1-i}
\left(w-\frac{\pi}{t},u\right) \d u = O\left(\frac{w^{\alpha}}{t^{k+1-
r}}\right).
\end{equation*}\vspace{.2pc}
\end{lem}

\begin{proof}
By an application of the mean value theorem for some $t+(1/w) < \xi <
\pi$,
\begin{align*}
&\int_{t+(1/w)}^{\pi} u^{r-i}(u-t)^{k-\alpha}S^{k,k+1-i}\left(w -
\frac{\pi}{t},u\right)\d u\\
&\quad\ = w^{\alpha-k}\int_{t+(1/w)}^{\xi} u^{r-i} S^{k,k+1-i}\left(w-
\frac{\pi}{t},u\right)\d u\\
&\quad\ = w^{\alpha-k}\left[u^{r-i}\ S^{k,k-i}\left(w-\frac{\pi}{t},
u\right)\right]_{u=t + (1/w)}^{\xi}\\
&\qquad\ \ -(r-i)w^{\alpha-k} \int_{t+(1/w)}^{\xi} u^{r-i-1}S^{k,k-i} \left(w-
\frac{\pi}{t},u\right)\d u\\
&\quad\ = w^{\alpha-k}
O\left\lbrace\frac{w^{k}}{\left(t+\frac{1}{w}\right)^{k+1-
r}}\right\rbrace + w^{\alpha-k} \int_{t+(1/w)}^{\xi} u^{r-i-1}
O\left(\frac{w^{k}}{u^{k+1-i}}\right)\d u\\
&\hskip 8.5cm\hbox{by Lemma~5(i)}\\[.2pc]
&\quad\ = O\left(\frac{w^{\alpha}}{t^{k+1-r}}\right) +
O\left(w^{\alpha}\int_{t+(1/w)}^{\xi} \frac{1}{u^{k+2-r}}\d
u\right)\\[.2pc]
&\quad\ = O \left(\frac{w^{\alpha}}{t^{k+1-r}}\right).
\end{align*}
\end{proof}

$\left.\right.$\vspace{-3.9pc}

\begin{lem}
For $i = 0, 1, 2, \ldots, m${\rm ,}
\begin{align*}
&\int_{t+(1/w)}^{\pi} u^{r-i} (u-t)^{k-\alpha}\ \d u \int_{1}^{w-
(\pi/t)} g_{i} (x,w,u)\d x\\[.2pc]
&\quad\ = O\left\lbrace\frac{w^{\alpha-
k}q^{k}\left(1-\frac{\pi}{wt}\right)}{t^{k+1-r}}\right\rbrace.
\end{align*}
\end{lem}

\begin{proof}
For some $1<\eta <w - (\pi/t)$, by an application of the mean value
theorem
\begin{align*}
&\int_{t+(1/w)}^{\pi} u^{r-i}(u-t)^{k-\alpha} \d u
\int_{t}^{w-(\pi/t)} g_{i} (x,w,u)\d x\\[.2pc]
&\quad\ = \frac{1}{k!} \int_{t-(1/w)}^{\pi} u^{r-i} (u-t)^{k-\alpha} \d
u \left[(-1)^{k}\left(\frac{\d}{\d x}\right)^{k}\ 
q_{\alpha}\left(\frac{x}{w}\right)\right]_{x=w-
(\pi/t)}\\[.2pc]
&\qquad\ \ \times \int_{\eta}^{w-(\pi/t)} \frac{\d}{\d x} S^{k,k+1-i}
(x,u)\d x\\[.2pc]
&\quad\ = \frac{q^{k}\left(1-\frac{\pi}{wt}\right)}{k!\ w^{k}}
\left\lbrace\int_{t+(1/w)}^{\pi} u^{r-i}(u-t)^{k-\alpha}S^{k,k+1-
i} \left(w-\frac{\pi}{t},u\right)\d u\right.\\[.2pc]
&\qquad\ \ \left. -\int_{t+(1/w)}^{\pi} u^{r-i} (u-t)^{k-\alpha}\ S^{k,k+1-
i}(\eta, u)\d u\right\rbrace\\[.2pc]
&\quad\ = O\left(\frac{q^{k}\left(1-\frac{\pi}{wt}\right)w^{\alpha}}{w^{k}\
t^{k+1-r}}\right)\\[.2pc]
&\quad\ = O \left(\frac{q^{k}\left(1-\frac{\pi}{wt}\right)w^{\alpha-k}}{t^{k+1-
r}}\right),
\end{align*}
since by Lemma~14, the first integral is $\displaystyle
O(w^{\alpha}/(t^{k+1-r}))$ and the second integral is dominated by
the first integral.
\end{proof}

\begin{lem}
For $i = 0,1,2,\ldots,m$ and $wu > \pi${\rm ,}
\begin{equation*}
\int_{w-(\pi/t)}^{w} g_{i} (x,w,u)\d x = O\left(w^{2}\ u^{-k+i}Q_{k}
\left(\frac{\pi}{wt}\right)\right).
\end{equation*}
\end{lem}

\begin{proof}
For $0\leq 1$, by use of Lemma~5(i),
\begin{align*}
&\int_{w-(\pi/t)}^{w} g_{i}(x,w,u)\d x\\[.2pc]
&\quad\ = \int_{w-(\pi/t)}^{w} \frac{(-1)^{k}}{(k-1)!}\left(\frac{\d}{\d x}\right)^{k}\ 
q_{\alpha}\left(\frac{x}{w}\right)S^{k-1,k+1-i} (x,u)\d x\\[.2pc]
&\quad\ = \frac{1}{(k-1)!}\int_{w-(\pi/t)}^{w} (-
1)^{k}\left(\frac{\d}{\d x}\right)^{k} \
q_{\alpha}\left(\frac{x}{w}\right) O\left(\frac{x^{k+1-
i}}{u^{k}}\right)\d x
\end{align*}
\begin{align*}
&\quad\ = O\left(\frac{w^{2-i}}{u^{k}}\int_{1-(\pi/wt)}^{1}
q^{k}(\theta) \d \theta\right)\\[.2pc]
&\quad\ = O\left(w^{2}u^{-k+i}\ Q_{k} \left(\frac{\pi}{wt}\right)\right)\quad
\hbox{as}\quad wu > \pi.
\end{align*}

For $i \geq 2$, by use of Lemma~5(i)
\begin{align*}
&\int_{w-(\pi/t)}^{w} g_{i}(x,w,u)\d x\\[.2pc]
&\quad\ = \int_{w-
(\pi/t)}^{w} \frac{(-1)^{k}}{(k-1)!}\left(\frac{\d}{\d x}\right)^{k}\ 
q_{\alpha}\left(\frac{x}{w}\right)\ S^{k-1,k+1-i} (x,u)\d x\\[.2pc]
&\quad\ = \frac{1}{(k-1)!}\int_{w-(\pi/t)}^{w} (-
1)^{k}\left(\frac{\d}{\d x}\right)^{k} \ 
q_{\alpha}\left(\frac{x}{w}\right) O\left(\frac{x^{k-1}}
{u^{k+2-i}}\right)\d x\\[.2pc]
&\quad\ = O\left(u^{-k-2+i} \int_{1-(\pi/wt)}^{1}
q^{k}(\theta)\d \theta\right)\\[.2pc]
&\quad\ = O\left(w^{2}u^{-k+i}\ Q_{k}\left(\frac{\pi}{wt}\right)\right)\quad
\hbox{as}\quad wu > \pi.
\end{align*}
Hence 
\begin{equation*}
\int_{w-(\pi/t)}^{w} g_{i}(x,w,u) \d x = O\left(w^{2}\ u^{-k+i}\
Q_{k}\left(\frac{\pi}{wt}\right)\right).
\end{equation*}
\end{proof}

\begin{lem}
For $i = 0,1,2,\ldots,m$ and $wt > \pi${\rm ,}
\begin{equation*}
\int_{t+(1/w)}^{\pi} u^{r-i} (u-t)^{k-\alpha} \frac{\d}{\d x}\
S^{k,k+1-i}(x,u) \d u = O\left(\frac{w^{\alpha}}{t^{k-r}}\right).
\end{equation*}
\end{lem}

\begin{proof}
Let $i = 0$. By mean value theorem for some $t + (1/w) < \xi < \pi$,
\begin{align*}
&\int_{t + (1/w)}^{\pi} u^{r - i} (u - t)^{k - \alpha}
\frac{\d}{\d x} S^{k, k + 1 - i} (x,u)\, \d u\\[.2pc]
&\quad \ = k \int_{t + (1/w)}^{\pi} u^{r} (u - t)^{k - \alpha} S^{k - 1,
k + 1} (x, u) \,\d u\\[.2pc]
&\quad \ = k w^{\alpha - k} \int_{t + (1/w)}^{\xi} u^{r} S^{k - 1, k + 1}
(x, u) \, \d u\\[.2pc]
&\quad \ = k w^{\alpha - k} \big[u^{r} S^{k - 1, k}(x, u)\big]^{\xi}_{u = t +
(1/w)}\\[.2pc]
&\qquad \ \ - k rw^{\alpha - k} \int_{t + (1/w)}^{\xi} u^{r - 1}
S^{k - 1, k} (x, u)\,\d u
\end{align*}
\begin{align*}
&\quad \ = kw^{\alpha - k} O \left\lbrace\frac{x^{k}}{(t + (1/w))^{k -
r}}\right\rbrace\\[.2pc]
&\qquad\ \  - k rw^{\alpha - k} \big[ u^{r - 1} S^{k - 1, k - 1}
(x, u)\big]^{\xi}_{u = t + (1/w)}\\[.1pc]
&\qquad\ \ + kr(r - 1) w^{\alpha - k} \int_{t + (1/w)}^{\xi} u^{r -
2} S^{k - 1, k - 1} (x, u)\,\d u \quad \hbox{by Lemma~5(i)}\\[.1pc]
&\quad \ = O \left(\frac{w^{\alpha}}{t^{k - r}}\right) + kr w^{\alpha - k}
O\left\lbrace \frac{x^{k - 1}}{\left(t + \frac{1}{w} \right)^{k - r + 1}}\right\rbrace\\[.1pc]
&\qquad\ \ + kr (r - 1) w^{\alpha - k} \int_{t + (1/w)}^{\xi} u^{r
- 2} O\left(\frac{x^{k - 1}}{u^{k}}\right) \ \d u\\[.2pc]
&\quad \ = O\left(\frac{w^{\alpha}}{t^{k - r}}\right) + O \left(\frac{w^{\alpha
- 1}}{t^{k - r +1}}\right)\\[.2pc]
&\quad \ = O\left(\frac{w^{\alpha}}{t^{k - r}}\right) \quad \hbox{as} \quad wt >
\pi.
\end{align*}
For $i > 1$, using the technique used in the proof of Lemma~14, it can
be proved that
\begin{align*}
&\int_{t + (1/w)}^{\pi} u^{r - i} (u - t)^{k -\alpha} \frac{\d}{\d
x}\ S^{k, k + 1 - i} (x, u) \,\d u\\[.1pc]
&\quad\ = O \left(\frac{w^{\alpha - 1}}{t^{k
- r + 1}}\right)\\[.1pc]
&\quad\ = O\left(\frac{w^{\alpha}}{t^{k - r}}\right) \quad \hbox{as} \quad wt
>\pi.
\end{align*}
This completes the proof of Lemma~17.\vspace{.5pc}
\end{proof}

\begin{lem}
For $i = 0, 1, 2, \ldots\!\,, m${\rm ,}
\begin{equation*}
\int_{1}^{\pi/t} \frac{\hbox{\rm \d} w}{w^{2}} \left\vert
\int_{t}^{\pi} u^{r - i} (u - t)^{k - \alpha} G_{i} (w, u) \,\hbox{\rm
\d} u\right\vert = O\left(\frac{1}{t^{\alpha - r}}\right).
\end{equation*}\vspace{.5pc}
\end{lem}

\begin{proof}
\begin{align*}
&\int_{1}^{\pi/t} \frac{\d w}{w^{2}} \left\vert
\int_{t}^{\pi} u^{r - i} (u - t)^{k - \alpha} G_{i} (w, u)\,\d
u\right\vert\\[.1pc]
&\quad \ \leq \int_{1}^{\pi/t} \frac{\d w}{w^{2}} \left\vert
\int_{t}^{t + (1/w)} u^{r - i} (u - t)^{k - \alpha} G_{i} (w,
u)\,\d u\right\vert\\[.1pc]
&\qquad\ \ + \int_{1}^{\pi/t} \frac{\d w}{w^{2}} \left\vert
\int_{\,t + (1/w)}^{\pi} u^{r - i} (u - t)^{k -\alpha} G_{i} (w,
u)\,\d u\right\vert\\[.1pc]
&\quad \ = \int_{1}^{\pi/t} \frac{\d w}{w^{2}} O(w^{\alpha - r+ 1}),
\quad \hbox{by Lemmas~11 and 12}\\[.1pc]
&\quad \ = O\left(\frac{1}{t^{\alpha - r}}\right).
\end{align*}
\end{proof}

$\left.\right.$\vspace{-3.15pc}

\begin{lem}
For $i = 0, 1, 2,\ldots\!\,, m$ and $wt > \pi${\rm ,}
\begin{align*}
&\int_{t}^{\pi} u^{r - i} (u - t)^{k - \alpha} G_{i} (w, u)\,\d u\\[.2pc]
&\quad\ = O\left(\frac{w^{\alpha - k} q^{k} \left(1 -
\frac{\pi}{wt}\right)}{t^{k - r + 1}}\right) + O\left(\frac{w^{\alpha -
k + 1}Q_{k} \left(\frac{\pi}{wt}\right)}{t^{k - r}}\right).
\end{align*}\vspace{.5pc}
\end{lem}

\begin{proof}
Using Lemma~10,
\begin{align*}
&\int_{t}^{\pi} u^{r - i} (u - t)^{k - \alpha} G_{i} (w, u)\,\d
u\\[.2pc]
&\quad \ = \int_{t}^{\pi} u^{r-i} (u - t)^{k - \alpha} \ \d u
\int_{1}^{w} g_{i}(x, w, u)\,\d x\\[.2pc]
&\quad \ = \int_{t}^{t + (1/w)} u^{r - i} (u - t)^{k - \alpha } \ \d u
\int_{1}^{w} g_{i} (x, w, u) \,\d x\\[.2pc]
&\qquad\ \  + \int_{t + (1/w)}^{\pi} u^{r - i} (u - t)^{k - \alpha}
\ \d u \int_{1}^{w} g_{i} (x, w, u)\,\d x\\[.2pc]
&\quad \ = J_{1} + J_{2}, \quad \hbox{say}.
\end{align*}

Using Lemmas~13 and 16,
\begin{align*}
J_{1} &= \int_{t}^{t + (1/w)} u^{r - i} (u - t)^{k - \alpha}\,\d
u \int_{1}^{w - (\pi/t)} g_{i} (x, w, u)\,\d x\\[.2pc]
&\quad \ + \int_{t}^{t + (1/w)} u^{r - i} (u - t)^{k -
\alpha}\,\d u \int_{w - (\pi/t)}^{w} g_{i} (x, w, u) \,\d x\\[.2pc]
&= O \left(\frac{w^{\alpha - k} q^{k} \left(1 - \frac{\pi}{wt}\right)}{t^{k - r
+ 1}}\right)\\[.2pc]
&\quad \  + O\left( \int_{t}^{t + (1/w)} u^{r - k} (u - t)^{k
- \alpha} w^{2} Q_{k} \left(\frac{\pi}{wt}\right)\,\d u\right)\\[.2pc]
&= O\left(\frac{w^{\alpha - k} q^{k} \left(1 -
\frac{\pi}{wt}\right)}{t^{k - r + 1}}\right) + O\left(\frac{w^{\alpha -
k + 1} Q_{k} \left(\frac{\pi}{wt}\right)}{t^{k - r}}\right) \quad
\hbox{as}\quad k \geq r
\end{align*}
and
\begin{align*}
J_{2} &= \int_{t + (1/w)}^{\pi} u^{r - i} (u - t)^{k - \alpha} \
\d u \int_{1}^{w - (\pi/t)} g_{i} (x, w, u) \,\d x\\[.2pc]
&\quad \ + \int_{t + (1/w)}^{\pi} u^{r - i} (u - t)^{k - \alpha} \
\d u \int_{w - (\pi/t)}^{w} g_{i} (x, w, u)\, \d x\\[.2pc]
&= O \left( \frac{w^{\alpha - k } q^{k} \left(1 -
\frac{\pi}{wt}\right)}{t^{k - r + 1}}\right) + \int_{w - (\pi/t)}^{w} \ \d x\\[.2pc]
&\quad \ \times \int_{t +
(1/w)}^{\pi} u^{r - i} (u - t)^{k - \alpha} g_{i} (x, w, u)\,\d x
\quad \hbox{by Lemma~15,}
\end{align*}
\begin{align*}
&= O \left( \frac{w^{\alpha - k} q^{k} \left(1 -
\frac{\pi}{wt}\right)}{t^{k - r + 1}}\right) 
+ \frac{1}{k!}\int_{w - (\pi/t)}^{w} (-1)^{k}
\left(\frac{\d}{\d x}\right)^{k} q_{\alpha} \left(\frac{x}{w}\right) \
\d x\\[.2pc]
&\quad \ \times \int_{t + (1/w)}^{\pi} u^{r - i} (u - t)^{k - \alpha}
\frac{\d}{\d x} S^{k, k + 1 -i} (x, u)\,\d x\\[.2pc]
&= O \left(\frac{w^{\alpha - k} q^{k} \left(1 -
\frac{\pi}{wt}\right)}{t^{k - r + 1}}\right)\\[.2pc]
&\quad \ + \frac{1}{k!} \int_{w - (\pi/t)}^{w} (-1)^{k} \left(\frac{\d}{\d x}\right)^{k}
q_{\alpha} \left(\frac{x}{w}\right) O \left(\frac{w^{\alpha}}{t^{k -
r}}\right)\,\d x\quad \hbox{by Lemma~17}\\[.2pc]
&= O \left(\frac{w^{\alpha - k} q^{k} \left(1 -
\frac{\pi}{wt}\right)}{t^{k - r + 1}}\right) + O\left(\frac{w^{\alpha - k +
1}}{t^{k - r}} \int_{1 - (\pi/wt)}^{1} q^{k} (\theta)\,
\d\theta\right)\\[.2pc]
&= O \left(\frac{w^{\alpha - k} q^{k} \left(1 -
\frac{\pi}{wt}\right)}{t^{k - r + 1}}\right) + O \left(\frac{w^{\alpha -
k + 1} Q_{k} \left(\frac{\pi}{wt}\right)}{t^{k - r}}\right).
\end{align*}
This completes the proof of Lemma~19.
\end{proof}

\begin{lem}
For $i = 0, 1, 2, \ldots\!\,, m${\rm ,}
\begin{equation*}
\int_{\pi/t}^{\infty} \frac{\d w}{w^{2}} \left\vert
\int_{t}^{\pi} u^{r - i} (u - t)^{k - \alpha} G_{i} (w, u)\,\d u
\right\vert = O \left(\frac{1}{t^{\alpha - r}}\right).
\end{equation*}\vspace{.5pc}
\end{lem}

\begin{proof}
By the use of Lemma~19, we get
\begin{align*}
&\int_{\pi/t}^{\infty} \frac{\d w}{w^{2}} \left\vert
\int_{t}^{\pi} u^{r - i} (u - t)^{k - \alpha} G_{i} (w, u)\,\d
u\right\vert\\[.2pc]
&\quad \ = O \left(\int_{\pi/t}^{\infty} \frac{w^{\alpha - k - 2} q^{k}
\left( 1- \frac{\pi}{wt}\right)}{t^{k - r + 1}}\,\d w\right) + O \left(
\int_{\pi /t}^{\infty} \frac{w^{\alpha - k - 1} Q_{k}
\left(\frac{\pi}{wt}\right)}{t^{k - r}}\,\d w\right)\\[.2pc]
&\quad \ = O \left(\frac{1}{t^{\alpha - r}} \int_{0}^{1}
\frac{q^{k}(\theta)}{(1 - \theta)^{\alpha - k}} \,\d \theta\right) +
O\left(\frac{1}{t^{\alpha - r}} \int_{0}^{1}
\frac{Q_{k}(u)}{u^{\alpha - k + 1}}\,\d u\right)\\[.2pc]
&\quad \ = O \left(\frac{1}{t^{\alpha - r}}\right) \quad \hbox{by Lemma~4.}
\end{align*}
\end{proof}

\section{Proof of the theorem}

\noindent{\it Proof of Theorem} 1.\hskip .5pc We have for $r \geq 1$,
\begin{align*}
\left(\frac{\d}{\d x}\right)^{r} B_{n} (x) &= \frac{2}{\pi}
\int_{0}^{\pi} \frac{(-1)^{r}}{2} \big\lbrace f(x + u) - (-
1)^{r} f(x - w)\big\rbrace\\[.2pc]
&\quad\ \times \left(\frac{\d}{\d u}\right)^{r} \sin nu\,\d u
\end{align*}
\begin{align*}
&= (-1)^{r} \frac{2}{\pi} \int_{0}^{\pi} h(u) u^{r}
\left(\frac{\d}{\d u}\right)^{r} \sin\,nu \,\d u\\[.2pc]
&\quad\ + (-1)^{r} \frac{2}{\pi} \int_{0}^{\pi} \frac{1}{2}
\big\lbrace P(u)\! -\! (-1)^{r} P(-u)\big\rbrace \left(\frac{\d}{\d
u}\right)^{r} \sin\,nu \,\d u\\[.2pc]
&= \alpha_{n} + \beta_{n}, \quad \hbox{say}.
\end{align*}
For the proof of our theorem it is enough to show that
\begin{equation*}
\sum \alpha_{n} \in \big| N_{q_{\alpha}}\big|
\end{equation*}
and
\begin{equation*}
\sum \beta_{n} \in \big| N_{q_{\alpha}}\big|.
\end{equation*}
Now
\renewcommand{\theequation}{\thesection\arabic{equation}}
\setcounter{equation}{0}
\begin{align}
n\alpha_{n} &= (-1)^{r} \frac{2}{\pi} \int^{\pi}_{0} n h(u) u^{r}
\left(\frac{\d}{\d u}\right)^{r} \sin\,nu \,\d u\nonumber\\[.2pc]
&= (-1)^{r + 1}\frac{2}{\pi} \int_{0}^{\pi} h(u) u^{r}
\left(\frac{\d}{\d u}\right)^{r + 1} \cos\,nu \,\d u\nonumber\\[.2pc]
&= (-1)^{r + 1} \frac{2}{\pi} \left\lbrack \sum\limits_{j =1}^{k - r} (-
1)^{j - 1} H_{j} (u) \left(\frac{\d}{\d u}\right)^{j-1} \right.\nonumber\\[.2pc]
&\quad\ \times \left\lbrace u^{r} \left.\left(\frac{\d}{\d u}\right)^{r + 1} \cos\,nu
\right\rbrace\right\rbrack_{u = 0}^{\pi}\nonumber\\[.2pc]
&\quad \ + (-1)^{k+1} \frac{2}{\pi}\int_{0}^{\pi} H_{k - r} (u)
\left(\frac{\d}{\d u}\right)^{k  - r} \left\lbrace u^{r}
\left(\frac{\d}{\d u}\right)^{r + 1} \cos\,nu\right\rbrace \,\d u\nonumber\\[.2pc]
&= J_{1}(n) + J_{2}(n), \quad \hbox{say}.
\end{align}
Since for $j = 1, 2, \ldots\!\,, k - r, H_{j} (+0) = O$ it is clear that
$J_{1}(n)$ is the sum of the terms containing $(-1)^{n} n^{p}$, where
$p$ is even and $r + 1\leq p \leq k$.

By the use of Lemma~9, for $p = 1, 2, \ldots\,\!, k$,
\begin{equation*}
\int_{1}^{\infty} \frac{\d w}{w^{2}} \left\vert \sum\limits_{n
\leq w} n^{p} (-1)^{n} q_{\alpha} \left(\frac{n}{w}\right) \right\vert <
\infty.
\end{equation*}
Hence
\begin{equation}
\int_{1}^{\infty} \frac{\d w}{w^{2}}\left\vert \sum\limits_{n
\leq w} J_{1}(n) q_{\alpha} \left(\frac{n}{w}\right)\right\vert <
\infty.
\end{equation}
Now
\begin{align*}
J_{2} (n) &= (-1)^{k + 1} \frac{2}{\pi} \int_{0}^{\pi} H_{k - r}
(u) \left(\frac{\d}{\d u}\right)^{k - r} \left\lbrace u^{r}
\left(\frac{\d}{\d u}\right)^{r + 1} \cos\,nu\right\rbrace \,\d u\\[.2pc]
&= \frac{2(-1)^{k + 1}}{\pi \Gamma (k - \alpha + 1)}
\int_{0}^{\pi} \left(\frac{\d}{\d u}\right)^{k - r} \left\lbrace
u^{r} \left(\frac{\d}{\d u}\right)^{r + 1} \cos\,nu\right\rbrace \,\d u\\[.2pc]
&\quad \ \times \int_{0}^{u} (u - t)^{k - \alpha} \d H_{\beta} (t) \quad
\hbox{by Lemma~1 as} \ \beta = \alpha - r \ \hbox{and} \ [\alpha] = k\\[.2pc]
& = \frac{2(-1)^{k + 1}}{\pi \Gamma (k - \alpha + 1)}
\int_{0}^{\pi} \d H_{\beta} (t) \int_{t}^{\pi} (u - t)^{k
- \alpha} \left(\frac{\d}{\d u}\right)^{k - r}\\[.2pc]
&\quad\ \times \left\lbrace u^{r}
\left(\frac{\d}{\d u}\right)^{r + 1} \cos\,nu \right\rbrace \,\d u\\[.2pc]
&= \frac{2(-1)^{k + 1}}{\pi \Gamma (k - \alpha + 1)}
\int_{0}^{\pi}\ \d H_{\beta} (t) \int_{t}^{\pi} (u - t)^{k
- \alpha}\\[.2pc]
&\quad \ \times \left\lbrace \sum\limits_{i = 0}^{m}\left(\begin{array}{c}
k - r\\
i
\end{array}\right) \left(\frac{\d}{\d u}\right)^{i} u^{r}
\left(\frac{\d}{\d u}\right)^{k + 1 -i} \cos\,nu\right\rbrace \,\d u\\[.2pc]
&\hskip 6cm \hbox{where}\ m = \min (k - r, r)\\[.2pc]
&= \frac{2(-1)^{k + 1}}{\pi \Gamma (k - \alpha + 1)} \sum\limits_{i =
0}^{m} \left(\begin{array}{c}
k - r\\
i
\end{array}\right) \frac{r!}{(r - i)!} \int_{0}^{\pi} \,\d
H_{\beta}(t)\\[.2pc]
&\quad \ \times \int_{t}^{\pi} (u - t)^{k - \alpha} u^{r - i}
\left(\frac{\d}{\d u}\right)^{k + 1 - i} \cos\,nu \,\d u.
\end{align*}
By the use of Lemmas~20 and 18,
\begin{align}
&\int_{1}^{\infty}\frac{\d w}{w^{2}} \left\vert \sum\limits_{n
\leq w} J_{2} (n) q_{\alpha} \left(\frac{n}{w}\right)\right\vert\nonumber\\[.2pc]
&\quad \ \leq \frac{2}{\pi \Gamma(k - \alpha + 1)}\sum\limits_{i =
0}^{m} \left(\begin{array}{c}
k - r\\
i
\end{array}\right) \frac{r!}{(r - i)!} \int_{0}^{\pi} \left\vert
\d H_{\beta} (t)\right\vert \int_{1}^{\infty} \frac{\d
w}{w^{2}}\nonumber\\[.2pc]
&\qquad\ \ \times \left\vert \int_{t}^{\pi} u^{r - i} (u -
t)^{k - \alpha} G_{i} (w, u)\,\d u\right\vert\nonumber\\[.2pc]
&\quad\ = \frac{2}{\pi \Gamma (k - \alpha + 1)} \sum\limits_{i = 0}^{m}
\left(\begin{array}{c}
k - r\\
i
\end{array}\right) \frac{r!}{(r - i)!} \int_{0}^{\pi} \left\vert
\d H_{\beta} (t)\right\vert\nonumber\\[.2pc]
&\qquad\ \ \times \left\lbrace \int_{1}^{\pi/t} \frac{\d w}{w^{2}}
\left\vert \int_{t}^{\pi} u^{r - 1} (u - t)^{k - \alpha} G_{i}
(w, u)\,\d u\right\vert\right.\nonumber\\[.2pc]
&\qquad\ \ \left.+ \int_{\pi/t}^{\infty} \frac{\d w}{w^{2}}
\left\vert\int_{t}^{\pi} u^{r - i} (u - t)^{k - \alpha} G_{i} (w,
u) \,\d u\right\vert \right\rbrace\nonumber
\end{align}
\begin{align}
&\quad\ = \frac{2}{\pi\Gamma (k - \alpha + 1)} \sum\limits_{i = 0}^{m} \left(\begin{array}{c}
k - i\\
i
\end{array}\right) \frac{r!}{(r - i)!} \int_{0}^{\pi} \left\vert
\d H_{\beta} (t)\right\vert O \left(\frac{1}{t^{\alpha - r}}\right)\nonumber\\[.2pc]
&\hskip 7.5cm \hbox{by Lemmas~18 and 20}\nonumber\\[.2pc]
&\quad\ = O\left(\sum\limits_{i = 0}^{m} \left(\begin{array}{c}
k - r\\
i
\end{array}\right) \frac{r!}{(r - i)!} \int_{0}^{\pi}
\frac{\big\vert\d H_{\beta} (t)\big\vert}{t^{\beta}}\right) \quad
\hbox{as} \ \alpha - r = \beta\nonumber\\[.2pc]
&\quad\ = O (1).
\end{align}
From (4.1), (4.2) and (4.3) it is clear that
\begin{equation*}
\sum \alpha_{n} \in |N_{q_{\alpha}}|.
\end{equation*}
Let $r$ be an odd number, i.e. $r = 2p+1$, where $p = 0,1,2,\ldots$ . Then
\begin{align*}
\beta_{n} &= -\frac{2}{\pi}\int_{0}^{\pi}\frac{1}{2} \{P(u) + P(-
u)\}\ \left(\frac{\d}{\d u}\right)^{2p+1} \sin\,nu\,\d u\\[.2pc]
&= (-1)^{p+1}\frac{2}{\pi} n^{2p+1} \int_{0}^{\pi}
\left(\sum\limits_{j=0}^{p} \frac{\theta_{2j} u^{2j}}{(2j)!}\right)
\cos\,nu\,\d u\\[.2pc]
&= (-1)^{p+1} \frac{2}{\pi}n^{2p+1} \sum\limits_{j=0}^{p}
\frac{\theta_{2j}}{(2j)!} \int_{0}^{u} u^{2j} \cos\,nu\,\d u\\[.2pc]
&= (-1)^{p+1}\frac{2}{\pi}n^{2p+1} \sum\limits_{j=1}^{p}
\frac{\theta_{2j}}{(2j)!} (-1)^{n}\\[.2pc]
&\quad \ \times \left(\sum\limits_{\mu=1}^{j} (-
1)^{\mu +1}n^{-2\mu}\pi^{2j-2\mu+1}\frac{(2j)!}{(2j-2\mu)!}\right)\\[.2pc]
&= 2(-1)^{n} \sum\limits_{\mu=1}^{p} (-1)^{p+\mu}n^{2p-2\mu+1}
\sum\limits_{j=\mu}^{p} \frac{\theta_{2j}}{(2j-2\mu)!} \pi^{2j-2\mu}.
\end{align*}
Let $r$ be an even number, i.e. $r = 2p$, where $p = 1,2,\ldots$ . Then
\begin{align*}
\beta_{n} &= \frac{2}{\pi}\int_{0}^{\pi}\frac{1}{2} \{P(u) - P(-
u)\}\ \left(\frac{\d}{\d u}\right)^{2p} \sin\,nu\,\d u\\[.2pc]
&= (-1)^{p} \frac{2}{\pi}n^{2p} \sum\limits_{j=1}^{p}
\frac{\theta_{2j-1}}{(2j-1)!} \int_{0}^{\pi} u^{2j-1} \sin\,nu\,\d u\\[.2pc]
&= (-1)^{p}\frac{2}{\pi}n^{2p} \sum\limits_{j=1}^{p}
\frac{\theta_{2j-1}}{(2j-1)!} (-1)^{n}\\[.2pc]
&\quad \ \times \left(\sum\limits_{\mu=1}^{j} (-
1)^{\mu -1}n^{-2\mu+1}\pi^{2j-2\mu+1}\frac{(2j-1)!}{(2j-2\mu)!}\right)\\[.2pc]
&= 2(-1)^{n} \sum\limits_{\mu=1}^{p} (-1)^{p+\mu-1}n^{2p-2\mu+1}
\sum\limits_{j=\mu}^{p} \frac{\theta_{2j-1}}{(2j-2\mu)!} \pi^{2j-2\mu}.
\end{align*}
So by the use of Lemma~9,
\begin{equation*}
\int_{1}^{\infty} \frac{\d w}{w^{2}} \left\vert \sum\limits_{n\leq w}
n\beta_{n} q_{\alpha} \left(\frac{n}{w}\right)\right\vert < \infty,
\end{equation*}
i.e. $\sum\beta_{n} \in \left\vert N_{q_{\alpha}}\right\vert$. This
terminates the proof of Theorem~1.

\end{document}